\documentclass[a4paper,oneside,10pt]{article}
\usepackage[english]{babel}
\usepackage[latin1]{inputenc}
\usepackage[T1]{fontenc}
\usepackage{amsmath}
\usepackage{amsfonts}
\usepackage{dsfont}
\usepackage{amssymb}
\title{Lelong-Skoda transform for compact Kähler manifolds and self-intersection inequalities}
\usepackage{dsfont}
\author{Gabriel Vigny}
\begin{document}
\newtheorem{theorem}{Theorem}[section]
\newtheorem{proposition}[theorem]{Proposition}
\newtheorem{defi}[theorem]{Definition}
\newtheorem{corollaire}[theorem]{Corollary}
\newtheorem{lemme}[theorem]{Lemma}
\newtheorem{Remark}[theorem]{Remark}
\maketitle

\begin{abstract}
Let $X$ be a compact Kähler manifold of dimension $k$ and $T$ be a positive closed current on $X$ of bidimension $(p,p)$ ($1\leq p < k-1$). We construct a continuous linear transform $\mathcal{L}_p(T)$ of $T$ which is a positive closed current on $X$ of bidimension $(k-1,k-1)$ which has the same Lelong numbers as $T$. We deduce from that construction self-intersection inequalities for positive closed currents of any bidegree. 
\end{abstract}


\section{Introduction}
Let $\mathcal{C}$ be a singular algebraic curve of degree $d$ in $\mathbb{P}^2$. Bézout's theorem implies that the number of intersections near all the singularities between $\mathcal{C}$ and a small perturbation of $\mathcal{C}$ using an automorphism satisfies:
\begin{eqnarray*}
\sum_{n}d_n(d_n-1) \leq d^2,
\end{eqnarray*}
where the $d_n$ are the multiplicity of the singularities (we refer the reader to the survey \cite{dem4} for more on that topic). Let $(X,\omega)$ be a compact Kähler manifold of dimension $k$. Here $\omega$ is a fixed Kähler form. We want to extend that inequality to the case where $\mathcal{C}$ is replaced by a variety of any bidimension. Two main difficulties appear: one cannot perturb the variety, and the self-intersection is not defined for varieties of small dimension. That is why we work on the more general case of positive closed currents on $X$. In the case of currents, the multiplicity of the singularities are replaced by the Lelong numbers. \\

Let $T$ be a positive closed current of bidimension $(p,p)$ in $X$. Let $E_c$ ($c>0$) denotes the set of point $z\in X$ where the Lelong number $\nu(T,z)$ of $T$ at $z$ is strictly larger than $c$. Siu's theorem implies that $E_c$ is an analytic subset of $X$ (possibly empty) of dimension $\leq p$ \cite{Siu}. Define $b_q:=\inf\{c>0, \ \text{dim} E_c \leq q \}$ and $b_{-1}:=\max_{x\in X} \nu(T,x)$. Then $0= b_{p}\leq \dots \leq b_{-1}$. For $c\in ]b_q,b_{q-1}]$, the dimension of $E_c$ is equal to $q$. So at $b_q$, there is a jump in the dimension for the analytic set $E_c$. Let $(Z_{q,n})_{n\geq 1}$ be the at most countable family of irreductible components of dimension $q$ of the $E_c$ for $c\in ]b_q,b_{q-1}]$; note that $(Z_{q,n})_{n\geq 1}$ is finite if $b_q>0$. Let $\nu_{q,n}:=\min_{x\in Z_{q,n}} \nu(T,x)\in ]b_q,b_{q-1}]$ be the generic Lelong number of $T$ on $Z_{q,n}$. Then we obtain:
 \begin{theorem}\label{1}
 Let $T$ be a	positive closed current of bidimension $(p,p)$ on $X$. Then with the above notation, for $q\leq p$ we have:
 \begin{eqnarray}\label{auto-intersection} 
 \sum_{n \geq 1} (\nu_{q,n}-b_{p})\dots(\nu_{q,n}-b_{q}) \|Z_{q,n}\| \leq C \|T\|^{p+1-q},
 \end{eqnarray}
 where $C$ is a constant which depends only on $(X,\omega)$, $\|T\|$ denotes the mass of the current $T$ and $\|Z_{q,n}\|$ denotes the mass of the current of integration on $Z_{q,n}$. 
 \end{theorem}
 Of course we can apply that inequality to the case where $T$ is the current of integration on subvarieties of $X$. Note that $\|Z_{q,n}\|$ is equal to $q!$ times the volume of $Z_{q,n}$. We refer to \cite{DS7} for a recent application of such inequalities in complex dynamics. In that paper, the authors prove that the normalized pull-back by $f^n$ of a generic hypersurface (in the Zariski sense) converges to the Green current of $f$ for a holomorphic endomorphism $f$ of $\mathbb{P}^k$. \\

For $p=k-1$, the inequality was proved  by Demailly using a regularization of quasi plurisubharmonic functions in \cite{dem1} (with $C=1$  when $X=\mathbb{P}^k$). Méo extended the result to the case of any bidimension for $X$ projective \cite{Meo}. The key point of his proof is that for a current of bidimension $(p,p)$ in $\mathbb{P}^k$, he constructed a positive closed current of bidegree $(1,1)$ having the same degree than $T$ and the same Lelong numbers. Such currents were constructed in $\mathbb{C}^k$ by Lelong and Skoda in \cite{Lel} and \cite{Sko}. They used a kernel in order to express the potential of the bidegree $(1,1)$ current. Méo used a geometric approach taking advantage of the fact that the family  of projective subspaces of dimension $p$ in $\mathbb{P}^{k}$ is very rich. The case of projective manifolds is then deduced from an embedding of $X$ in some $\mathbb{P}^N$. \\

Theorem \ref{1} generalizes this result to the Kähler case. For that, we use an approach inspired by \cite{DS5} (see also \cite{BGS}, \cite{GS} and \cite{DS8}). Let $T$ be a positive closed current of bidimension $(p,p)$ on $X$ for $1\leq p<k-1$. Consider the canonical projections $\pi_1, \ \pi_2:X\times X \to X$ and let $\pi:\widetilde{X\times X}\to X\times X$ be the blow-up of $X\times X$ along the diagonal $\Delta$ of $X\times X$. We modify the pull-back by $\pi_2$ of the current $T$ by multiplying it by a suitable form $\Theta$ and then we push it forward to $X$ by $\pi_1$. 
The form $\Theta$ is smooth outside $\Delta$. It is defined as $(\pi_*(\widetilde{\Omega}))^{p+1}$ where $\widetilde{\Omega}$ is a Kähler form in $\widetilde{X\times X}$. This is what we call the \emph{Lelong-Skoda transform} $\mathcal{L}_p(T)$ of the current $T$. The Lelong-Skoda transform is a continuous linear operator which sends positive closed currents of bidimension $(p,p)$ to positive closed currents of bidegree $(1,1)$. The Lelong-Skoda transform is linear in the sense that it is linear on the space of currents spanned by positive closed currents. We show in Theorem \ref{main} below that $\mathcal{L}_p(T)$ has the same Lelong number than $T$ at every point. We then prove Theorem \ref{1} in the same way than in \cite{Meo}. We will also give some properties of the Lelong-Skoda transform. Finally, we extend the results to the case of harmonic currents.

\section{Lelong-Skoda transform}\label{geometric}
For $z=(z_1,\dots,z_m)\in \mathbb{C}^m$, we use the notation $|z|:=(\sum_j |z_j|^2)^{1/2}$.
Consider a compact Kähler manifold $Y$ of dimension $m$ endowed with a Kähler form $\theta$. We denote  by $\mathcal{C}_p(Y)$ the cone of positive closed currents of bidimension $(p,p)$ in $Y$ on which we consider the usual weak topology. For $T\in \mathcal{C}_p(Y)$, let $\|T\|:=\int_Y T \wedge \theta^{p}$ be the \emph{mass} of $T$. Note that this mass depends only on the class of $T$ in the Hodge cohomology group $H^{m-p,m-p}(Y,\mathbb{C})$. We refer the reader to the book of Voisin \cite{voi} for basics on compact Kähler manifolds.\\

 We recall some basic facts on Lelong numbers. See \cite{dem3} for proofs and details and also \cite{FS2} for results on intersection theory of currents. For $a \in Y$ and $T\in \mathcal{C}_p(Y)$, we consider a chart $V$ around $a$ in which the coordinates are given by $x$. For simplicity, assume $a$ is given by $0$ here. Then for $r>0$, the positive measure $T(x)\wedge (dd^c\log |x|)^p$ is well defined by \cite{dem3} and we define the quantities:
\begin{eqnarray*}
\nu(T,0,r) &:=& \int_{|x|<r} T(x)\wedge (dd^c\log |x|)^p\\
\nu(T,0)   &:=& \lim_{r\to 0} \int_{|x|<r} T(x)\wedge (dd^c\log |x|)^p.
\end{eqnarray*}
The quantity $\nu(T,0)$ is called the \emph{Lelong number} of $T$ at $0$ and does not depend on the choice of the coordinates $x$ (so it is well defined on manifolds). It is the mass at $0$ of the measure $T(x)\wedge (dd^c\log |x|)^p $. We have the equivalent definition:
$$ \nu(T,0)= \lim_{r\to 0} r^{-2p} \int_{|x|<r} T \wedge (\frac{1}{2}dd^c |x|^2)^p.$$
It follows that the function $S\to \nu(S,0)$ is a linear form from the vector space spanned by positive closed currents, so the Lelong number is well defined for a difference of positive closed currents.

We give now a geometrical interpretation of the Lelong number that we will use latter. Let $\imath:\widetilde{Y} \to Y$ denote the blow-up of $Y$ at $a$ and $H:=\imath^{-1}(a)$ the exceptional divisor. Then $H \simeq \mathbb{P}^{m-1}$ and we put on $H$ the Fubini-Study form $\omega_{FS}$ normalized by $\int_{\mathbb{P}^{m-1}} \omega_{FS}^{m-1}=1$. Let $V$ and $x$ be as above. Consider  a sequence $(T_n)$ of smooth positive closed currents of bidimension $(p,p)$ in $V$ with $T_n \to T$ in the sense of currents in $V$. We can obtain $(T_n)$ using a convolution operator. Let $\widetilde{T}$ be a cluster value of the (bounded) sequence $(\imath^*(T_n))$. In bidegree $(1,1)$, $\widetilde{T}$ is unique and does not depend on $(T_n)$, this is not true in higher bidegree. Then we have the characterization:
\begin{lemme}\label{interpretation}
With the above notation, the Lelong number $\nu(T,0)$ is the mass of $\widetilde{T}$ on $H$. In particular, that mass does not depend on the choice of $(T_n)$.
\end{lemme}
\emph{Proof.} Observe that $\imath^{-1}(V)=\{(x,[u])\in V\times \mathbb{P}^{m-1}, \quad x \in [u] \}$. Write $x=(x_1,\dots,x_m)$ and $u=[u_1: \dots: u_m]$. Let $\jmath:\imath^{-1}(V)\to V\times \mathbb{P}^{m-1}$ be the canonical holomorphic injection. Pull-backs by $\jmath$ of positive closed smooth forms are positive closed smooth forms. Let $p:V\times \mathbb{P}^{m-1} \to \mathbb{P}^{m-1}$ be the projection on the second factor and recall that $\omega_{FS}$ is the  Fubini-Study form on $\mathbb{P}^{m-1}$. We endow $V\times \mathbb{P}^{m-1}$ with the form $p^*(\omega_{FS})$. We consider the form $\jmath^*(p^*(\omega_{FS}))$ on $\imath^{-1}(V)$. The restriction of $\jmath^*(p^*(\omega_{FS}))$ to $H$ is indeed $\omega_{FS}$ so we write for simplicity $\omega_{FS}([u])$ instead of $\jmath^*(p^*(\omega_{FS}))(x,[u])$, and we drop the $[u]$ when there can be no confusion.
In the chart of $\imath^{-1}(V)$ where $u_1=1$, write $u=[1:u_2:\dots:u_m]$ and consider the coordinates $(X_1,\dots, X_m)$ given by $X_1=x_1, \ X_2=u_2, \dots, \ X_m=u_m$ (so $X_i=x_i/x_1$ for $i>1$) in which $H$ is given by $(X_1=0)$. In these coordinates, the form $\omega_{FS}$ is:
$$\omega_{FS}([u])= \frac{1}{2} dd^c \log(1+|X_2|^2+\dots+|X_m|^2). $$
Since $dd^c \log |x_1|=0$ for $x_1\neq 0$ and $(dd^c\log |x|)^p$ does not charge the point $x=0$ for $p<k$, we have that:
\begin{eqnarray*}
\imath_*(\omega_{FS}([u]))&=&dd^c\log |x|\\
\imath_*(\omega_{FS}([u])^p)&=&(dd^c\log |x|)^p.
\end{eqnarray*}
Observe that we consider the push-forward of a smooth form and not $\imath^*(dd^c \log |x|)$ which is the pull-back of a current; indeed to have the "good" continuity's properties one should write $\imath^*(dd^c \log |x|)=\omega_{FS}([u])+[H]$. Consider a smooth psh function $f$ that coincides with $\log |x|$ for $|x|>r/2 $. Then Stokes formula gives that $\nu(T,0,r')=\|T\wedge (dd^c f)^p\|_{B_{r'}}$ for $r'>r/2$ ($B_{r'}$ is the ball of center $0$ and radius $r'$).  Then $T_n\wedge (dd^c f)^p\to T\wedge (dd^c f)^p$ in the sense of measures. In particular, for $\varepsilon>0$ small enough, we have:
$$\nu(T,0,(1-2\varepsilon)r)\leq \lim_{n\to \infty} \int_{B_{(1-\varepsilon)r}} T_n\wedge (dd^c f)^p\leq \nu(T,0,r).$$
Again, Stokes formula gives that $\int_{B_{(1-\varepsilon)r}} T_n\wedge (dd^c f)^p= \nu(T_n,0,(1-\varepsilon)r)$, so:
\begin{eqnarray}\label{trick}
\nu(T,0)= \lim_{r\to 0} \lim_{n \to \infty} \nu(T_n,0,r). 
\end{eqnarray}
Now, for $T_n$ smooth the current $\imath^*(T_n)$ is a well defined smooth form so:
\begin{eqnarray*}
\nu(T_n,0,r)&=&\int_{|x|<r} T_n(x)\wedge (dd^c\log |x|)^p \\
            &=&\int_{\imath^{-1}(|x|<r)} \imath^*(T_n(x)) \wedge \omega_{FS}^p.
\end{eqnarray*}
By definition of weak convergence we have for a set of $r$ of full measure that:
$$\lim_{n\to \infty} \left(\int_{\imath^{-1}(|x|<r)} \imath^*(T_n(x)) \wedge \omega_{FS}^p\right)\to \int_{\imath^{-1}(|x|<r)} \widetilde{T} \wedge \omega_{FS}^p.$$
So:
\begin{eqnarray*}
\nu(T,y)= \lim_{r\to 0} \int_{\imath^{-1}(|x|<r)} \widetilde{T} \wedge \omega_{FS}^p.
\end{eqnarray*}
The restriction to $H$ of $\omega_{FS}^p$ is the Fubini-Study form on $H$ at the power $p$. We decompose $\widetilde{T}= \widetilde{T}_1 +\widetilde{T}_2$ where $\widetilde{T}_1$ is the restriction of $\widetilde{T}$ to $H$ and $\widetilde{T}_2$ the restriction of $\widetilde{T}$ to $i^{-1}(V)\backslash H$ ($\widetilde{T}_2$ can also be defined as the trivial extension of $\imath^*(T)$ defined outside $H$). In other words:
\begin{eqnarray*}
\nu(T,0)=  \langle \widetilde{T}_1 , \omega_{FS}^p \rangle,
\end{eqnarray*}  
so $\nu(T,0)$ can be interpreted geometrically as the mass of the current $\widetilde{T}$ on the exceptional divisor. $\Box$ \\

 Let $X$ be a compact Kähler manifold of dimension $k$. Let $\Delta$ denote the diagonal in $X\times X$ and $\pi:\widetilde{X\times X} \to X\times X$ the blow-up of $X\times X$ along $\Delta$. Then $\widetilde{X\times X}$ is a compact Kähler manifold by Blanchard \cite{Bla}. Fix a Kähler form $\widetilde{\Omega}$ on $\widetilde{X\times X}$ and define $\Omega:=\pi_* (\widetilde{\Omega})$. Then $\Omega$ is a positive closed $(1,1)$ current smooth outside $\Delta$ and that does not charge $\Delta$. Observe that for $p < k-1$, $\Omega^{p+1}=\pi_* (\widetilde{\Omega}^{p+1})$ since those two positive closed currents coincide outside $\Delta$ and cannot charge $\Delta$ since they are of bidimension $>k=\text{dim}(\Delta)$. Define $\widetilde{\Delta}:=\pi^{-1}(\Delta)$.

Let $\pi_1$ and $\pi_2$ denote the canonical projections from $X\times X$ to its factors and denote $\widetilde{\pi_i}:=\pi_i\circ \pi$. Then $\pi_i$, $\widetilde{\pi_i}$ are submersions (see \cite{DS5}) so push-forward and pull-back by $\pi_i$ and $\widetilde{\pi_i}$ of positive closed currents are well defined  operators which are continuous for the topology of currents. The mass of $\widetilde{\Omega}\wedge [\widetilde{\pi_2}^{-1}(y)]$ does not depend on $y\in X$ for cohomological reasons and is positive since $\widetilde{\Omega}$ is a Kähler form. So renormalizing $\widetilde{\Omega}$ if necessary, we can assume that this mass is equal to $1$. The following result gives a description of the singularities of $\Omega$ near $\Delta$ and explains our choice of $\Omega$ (see also the definition of Lelong number):
\begin{proposition}\label{equivalent}
The form $\Omega$ admits locally a potential $\varphi$ that is $dd^c \varphi=\Omega$ where $\varphi$ is a psh function such that $\varphi(*)-\log {\rm dist}(*,\Delta)$ is bounded.
\end{proposition}
The problem is local so we let $U$ be a chart of $X\times X$ in which we consider the coordinates $(z,w)$ such that $\Delta$ is given by $z=0$ here and $\pi_2(z,w)=w$. So $\text{dist}(*,\Delta)$ is equivalent to $|z|$ here. Then $\pi^{-1}(U)=\{ (z,w,[u])\in U \times \mathbb{P}^{k-1}, \ z \in [u] \}$ and $\widetilde{\Delta}$ is given by $(z=0)$ here. 

As above, consider on $\pi^{-1}(U)$ the smooth form $\omega_{FS}([u])$ then $\widetilde{\alpha}:=\widetilde{\Omega}-\omega_{FS}([u])$ is a smooth closed form of bidegree $(1,1)$ such that $\int_{\mathbb{P}^{k-1}_w} \widetilde{\alpha}\wedge \widetilde{\Omega}^{k-2}=0$ (where $\mathbb{P}^{k-1}_w:=\{(0,w,[u]), \  [u]\in \mathbb{P}^{k-1}\}$). The idea is to show that $\widetilde{\alpha}$ is exact so it can be written $dd^c g$ for $g$ a continuous function. We will need the following lemma:
\begin{lemme}
Any closed form $\widetilde{\alpha}$ of bidegree $(1,1)$ on $\pi^{-1}(U)$  such that $\int_{\mathbb{P}^{k-1}_w} \widetilde{\alpha} \wedge \omega_{FS}^{k-2}=0$ for all $w$ is in fact exact.
\end{lemme}  
\emph{Proof.} We can assume that $\widetilde{\alpha}$ is real. Using Poincaré's lemma in the non-compact case, the only obstruction would be that there exists a compact oriented surface $S$ and a smooth function $f:S\to \pi^{-1}(U)$ such that $\int_{S} f^*(\widetilde{\alpha}) \neq 0$. To simplify the notations, assume that $(0,0)\in U$. If $i:\pi^{-1}(U)\times [0,1]\to \pi^{-1}(U) $ and $j:\pi^{-1}(U)\to \pi^{-1}(U) $ are the maps defined by $i((z:w:[u]),t)=(tz,tw,[u])$ and $j(z:w:[u])=(0,0,[u])$, we obtain by homotopy that $\int_{S} (j\circ f)^*(\widetilde{\alpha}) \neq 0$ which contradicts our hypothesis.  $\Box$ \\

\noindent\emph{End of the proof of Proposition \ref{equivalent}.} Let $p_w: \mathbb{P}^{k-1}_w \to \pi^{-1}(U)$ be the canonical injection. We denote by $\widetilde{\alpha}_w$ the form  $p_w^*(\widetilde{\alpha})$ on $\mathbb{P}^{k-1}_w$. Observe that $\int_{\mathbb{P}^{k-1}_w} \widetilde{\alpha} \wedge \omega_{FS}^{k-2}=\int_{\mathbb{P}^{k-1}_w} \widetilde{\alpha}_w \wedge \omega_{FS}^{k-2}$.

Then $\widetilde{\alpha}_w$ is a smooth exact form, hence it is $dd^c$ exact. So there exists a real smooth function $v_w$ such that $dd^c v_w=\widetilde{\alpha}_w$ (observe that any solution of $dd^c v_w=\widetilde{\alpha}_w$ is smooth). Since $v_w$ can be written via a kernel with coefficients in $L^1$ smooth outside the diagonal of $(\mathbb{P}^{k-1}_w)^2$, we deduce that $v(z,w,[u])=v_w([u])$ is continuous in all variables and smooth in the $w$ variable. A more precise analysis of the singularities of the kernel shows that its gradient has also coefficients in $L^1$ so $v(z,w,[u])$ is in $\mathcal{C}^1$ (see the proof of Proposition 2.3.1 in \cite{DS6}, we believe that we can choose $v$ smooth but we do not know how to prove it). Still the current $dd^cv$ is represented by a continuous form (we use the integral representation of $v$ and the fact that $\widetilde{\alpha}_w$ is smooth hence continuous in both variables). Replacing $\widetilde{\alpha}$ by $\widetilde{\Omega}-\omega_{FS}([u])-dd^c v$ in $\pi^{-1}(U)$, we can assume that $\widetilde{\alpha}_w=0$ for all $w \in \pi_2(U)$. 

Using the previous lemma, we can write $\widetilde{\alpha}=d \widetilde{\beta}$ where $\widetilde{\beta}$ is a  form in $\mathcal{C}^1$. We decompose $\widetilde{\beta}$ as $\widetilde{\beta}_{1,0}+\widetilde{\beta}_{0,1}$, its $(1,0)$  and $(0,1)$ components. We deduce from the equation $\widetilde{\alpha}=d\widetilde{\beta}$ that $\widetilde{\beta}_{0,1}$ is $\bar{\partial}$-closed and since  $\widetilde{\alpha}$ is real $\widetilde{\beta}_{0,1}=\overline{\widetilde{\beta}_{1,0}}$.  The form $p^*_w(\widetilde{\beta}_{0,1})$ is $\bar{\partial}$-closed and $\partial$-closed since $2\partial \widetilde{\beta}_{0,1}=\widetilde{\alpha}$ (we use here that $\widetilde{\alpha}_w=0$). So  $p^*_w(\widetilde{\beta}_{0,1})$ is a closed $(0,1)$-form on $\mathbb{P}^{k-1}$: it is equal to zero.

We work in the chart where $|u_1|>1/2 \max_i |u_i|$  so we take $u_1=1$. We consider the coordinates given by $(Z_1,\dots,Z_n,w)$ where:  
$$Z_1=z_1, \quad Z_2= z_2/z_1, \ \dots, \quad Z_k= z_k/z_1.$$
The form $\widetilde{\beta}_{0,1}$ can be written as $\sum_i \widetilde{\beta_i}d\overline{Z_i}+ \sum_i \widetilde{\beta'_i}d\overline{w_i}$ and since $p^*_w(\widetilde{\beta}_{0,1})=0$, we have that $\widetilde{\beta_i}(0,Z_2,\dots,Z_k,w)=0$ for all $i\geq 2$ while the other coefficients are bounded. Since the coefficients are in $\mathcal{C}^1$, we have that $|\widetilde{\beta_i}|\leq C |Z_1|$ for $i\geq 2$. Let $\beta_{0,1}=\pi_*(\widetilde{\beta}_{0,1})$, then in the local coordinates, it can be written:
$$\beta_{0,1}=\widetilde{\beta}_1dz_1+\sum_{i\geq 2} \widetilde{\beta_i}d(\frac{z_i}{z_1})+ \sum_i \widetilde{\beta'_i}dw_i.$$
We write $d(z_i/z_1)=dz_i/z_1-z_i/z^2_1dz_1$. Using the fact that the coefficients  $\widetilde{\beta_i}$ are bounded by  $C|z_1|$ and $|z_i|<2|z_1|$ we get that $\beta_{0,1}$ has coefficients in $L^\infty$. In particular, $\beta_{0,1}$  is a well defined $(0,1)$ current $\bar{\partial}$-closed with $L^\infty$-coefficients. So taking $U$ strictly pseudoconvex, we can solve $\bar{\partial} u= \beta_{0,1}$ with $u$ a continuous function \cite{HL1} ($L^\infty$ is enough for our purpose). Let $g:=-i\pi(u-\bar{u})$ it is a real bounded function and it satisfies the following identities for $\alpha:=\pi_*(\widetilde{\alpha})$ and $\beta:=\pi_*(\widetilde{\alpha})$:
$$dd^c g= \partial \bar{\partial} (u-\bar{u})= \partial \beta_{0,1}+\overline{\partial \beta_{1,0}}= d \beta=\alpha= \Omega-dd^c\log|z|-dd^c \pi_*(v).$$
This implies that $\Omega=dd^c(\log |z|+g-\pi_*(v))$ and gives the proposition ($\pi_*(v)$ is only bounded). \nolinebreak $\Box$\\

 Let $T \in \mathcal{C}_p(X)$ with $1\leq p< k-1$ (there is nothing to be done for $p=k-1$). The current $(\widetilde{\pi_2})^*(T) \wedge \widetilde{\Omega}^{p+1}$ is a well defined element of $\mathcal{C}_{k-1}(\widetilde{X\times X})$. It is of finite mass and coincides with $(\pi_2)^*(T) \wedge \Omega^{p+1}$ outside $\widetilde{\pi}^{-1}(\Delta)$. So $(\pi_2)^*(T) \wedge \Omega ^{p+1}$ is well defined on $X\times X$ as the trivial extension over $\Delta$ of the above current. Consequently, the current: 
 $$T_{LS}:= (\pi_1)_*((\pi_2)^*(T) \wedge \Omega ^{p+1})$$
belongs to $\mathcal{C}_{k-1}(X)$. That gives us the following notion:
\begin{defi}
For $p <k-1$, define the \emph{Lelong-Skoda transform} $\mathcal{L}_p$ from $\mathcal{C}_p(X)$ to $\mathcal{C}_{k-1}(X)$ by:
$$\mathcal{L}_p(T):=T_{LS}.$$
\end{defi}
Of course, the operator $\mathcal{L}_p$ depends on the choice of $\widetilde{\Omega}$.

\section{Properties of the Lelong-Skoda transform}\label{properties}
This section is devoted to prove the following result:
\begin{theorem}\label{main}
The Lelong-Skoda transform $\mathcal{L}_p$ is a continuous linear operator. In particular, there exists a constant $C>0$ such that $\|\mathcal{L}_p(T)\|\leq C \| T\|$. Moreover, it preserves the Lelong number, i.e. $\mathcal{L}_p(T)$ and $T$ have the same Lelong number at every point. 
\end{theorem} 
 Since the set $\mathcal{C}_p(X)$ is a convex cone, we should not speak of linear operator. Nevertheless, the transform  $\mathcal{L}_p$ can be extended to a linear operator on the space of currents spanned by positive closed currents. For a current $S$ which is the difference of positive closed currents, we can write:
 $$\| S\|:= \inf (\|T'\|+\|T''\|)$$ 
 where the infimum is taken over all the decompositions $S=T' - T''$ where $T'$ and $T''$ are positive closed currents. 
 We use the following lemma to prove the continuity:
\begin{lemme}\label{Fubini}
The transform $\mathcal{L}_p$ satisfies : $\mathcal{L}_p(T)=(\widetilde{\pi_1})_*((\widetilde{\pi_2})^*(T) \wedge \widetilde{\Omega}^{p+1})$. 
\end{lemme}
\emph{Proof.} Since $\widetilde{\Omega}$ is smooth, it is sufficient to show that $(\widetilde{\pi_2})^*(T)$ does not charge $\widetilde{\Delta}$. 
The problem is local so we let $U$ be a chart of $X\times X$ in which we consider the coordinates $(z,w)$ such that $\Delta$ is given by $z=0$ here and $\pi_2(z,w)=w$. Then $\pi^{-1}(U)=\{ (z,w,[u])\in U \times \mathbb{P}^{k-1}, \quad z \in [u] \}$ and $\widetilde{\Delta}$ is given by $z=0$ here, so $\pi^{-1}(U)$ is locally a product. So $\widetilde{\pi_2}^*$ of a current is just integration on fibers. We want to know if $(\widetilde{\pi_2})^*(T)$ charges $(z=0)$ which is impossible by Fubini's theorem. More precisely, let $B_r$ denote the ball of center $0$ and radius $r$ in $\mathbb{C}^{k}$. We can reduce ourselves to the case where $\pi^{-1}(U)$ is of the form $\widehat{B_{r_1}}\times B_{r_2}$ where $\widehat{B_{r_1}}$ is the blow-up of $B_{r_1}$ at $z=0$. Let $\omega_{\widehat{B_{r_1}}}$ and $\omega_{B_{r_2}}$ denote Kähler form on $\widehat{B_{r_1}}$ and $B_{r_2}$. Then:
\begin{eqnarray*}
\|(\widetilde{\pi_2})^*(T)\|_{\widehat{B_{r}}\times B_{r_2}}&=& \binom{k+p}{k}\int_{B_{r_2}}T(w)\wedge \omega^{p}_{B_{r_2}}  \int_{\widehat{B_{r}}} \omega^k_{\widehat{B_{r_1}}} \\
                                                            &=& \binom{k+p}{k}\|T\|_{B_{r_2}} \omega^k_{\widehat{B_{r_1}}}(\widehat{B_{r}}).
\end{eqnarray*}
 And the lemma follows from letting $r$ goes to $0$. $\Box$\\

\noindent The lemma implies that $\mathcal{L}_p$ is a continuous linear operator since $\widetilde{\Omega}^{p+1}$ is smooth and pull-back and push-forward by the submersions $\widetilde{\pi_1}$ and $\widetilde{\pi_2}$ are continuous operators on (positive closed) currents. \\

\noindent\emph{Proof of Theorem \ref{main}.} We want to interpret the mass of $\mu$ at $(0,0)$ as in Lemma \ref{interpretation}. For that purpose, we will pull-back some integrals on some suitable blow-ups of $X\times X$ in order to desingularize the forms $\Omega$, $dd^c \log |x|$ and $dd^c \log |(x,y)|$.\\

\noindent  Let $p_1:\widehat{X\times X}\to X\times X$ be the blow-up of $X\times X$ at $(0,0)$. Consider the system of local coordinates $(z,w)$ in the neighborhood $V\times V$ of $(0,0)$ in $X\times X$ given by $(z,w):=(x-y,x)$ for $(x,y)\in V\times V$. Then :
$$\widehat{V\times V}:=p_1^{-1}(V\times V)=\{(z,w, [u:v]) \in V\times V \times \mathbb{P}^{2k-1}, \ (z,w) \in [u:v] \}, $$  
where $z=(z_1,\dots,z_k)$, $w=(w_1,\dots,w_k)$, $[u:v]=[u_1:\dots:u_k:v_1:\dots:v_k]$. In $\widehat{V\times V}$, there exists a smooth form that we  denote by $\omega_{FS,2k-1}$ such that $dd^c \log |(z,w)|=(p_1)_*(\omega_{FS,2k-1})$ (that is what we mean by desingularization). But we cannot do that for $dd^c \log |w|$ nor for $\Omega$ so we need to blow-up once more. \\

\noindent  Consider the smooth submanifold $M$ of $\widehat{X\times X}$ given by $\{u=0\} \cup \{v=0 \}$ in $\widehat{V\times V}$ and by $p_1^{-1}(\Delta \cup (\{0\}\times X))$ outside $\widehat{V\times V}$. It is the disjoint union of two submanifolds which are the strict transform of $(z=0)$ and $(w=0)$. In other words, the blow-up $p_1$ desingularizes the analytic set $\Delta \cup \{x=0\}$. Let $p_2:\widehat{\widehat{X\times X}}\to \widehat{X\times X}$ denote the blow-up of $\widehat{X\times X}$ along $M$. So:
\begin{eqnarray*}
\widehat{\widehat{V\times V}}:=p_2^{-1}(\widehat{V\times V})= &\{ (z,w,[u:v],[u'],[v'])\in (V\times V) \times \mathbb{P}^{2k-1}\times \mathbb{P}^{k-1}\times \mathbb{P}^{k-1}, \\
                               &(z,w)\in [u:v], \ u \in [u'], \ v \in [v'] \}.
\end{eqnarray*}
Observe that $z\in [u']$. So we can define the holomorphic projection $P:\widehat{\widehat{X\times X}}\to \widetilde{X\times X}$  which is the identity outside $\widehat{\widehat{V\times V}}$ and given in $\widehat{\widehat{V\times V}}$ by:
$$P: (z,w,[u:v],[u'],[v']) \mapsto (z,w,[u']).$$
The form $P^*(\widetilde{\Omega})$ is a well defined positive smooth form that we will simply write $\widetilde{\Omega}$. \\

\noindent Define finally $\pi': \widetilde{X\times X}' \to X \times X$ the blow-up of $X \times X$ along $\{x=0\}$. Observe that $w\in [v']$. So we can define the holomorphic projection $P':\widehat{\widehat{X\times X}}\to \widetilde{X\times X}'$ which is given in $\widehat{\widehat{V\times V}}$ by:
$$P': (z,w,[u:v],[u'],[v']) \mapsto (z,w,[v']).$$
 Let $\omega_{FS,k-1}$ denote the smooth form on $(\pi')^{-1}(V\times X)$ which is the pull-pack of the Fubini-Study form on $\mathbb{P}^{k-1}$ to $(\pi')^{-1}(V\times X)$ (observe that there may not be a global mapping $\widetilde{X\times X}'\to X\times X \times \mathbb{P}^{k-1}$ so we cannot speak of $\omega_{FS,k-1}$ on $\widetilde{X\times X}'$). Then the form $(P')^*(\omega_{FS,k-1})$ is a well defined positive smooth form on $(p_1\circ p_2)^{-1}(V\times X)$ that we will simply write $\omega_{FS,k-1}$. \\

\noindent We define $E:=(p_1\circ p_2)^{-1}(\{(0,0) \})=\{((0,0),[u:v], [u'],[v']) \in \widehat{\widehat{V\times V}}\}$ the fiber of $(0,0)$ in $\widehat{\widehat{X\times X}}$. Note that $E$ can be considered as a blow-up of $\mathbb{P}^{2k-1}\simeq p_1^{-1}(0,0)$ along two disjoint subspaces of dimension $k-1$. Let $\Pi:E \to \mathbb{P}^{2k-1}$ be the canonical projection defined by:
$$\Pi: ([u:v],[u'],[v']) \mapsto [u:v].$$
So $\Pi$ is just the restriction of $p_2$ to $E$ if we identify $\mathbb{P}^{2k-1}$ to $p_1^{-1}(0,0)$. Finally let $\widehat{\pi_i}:=\pi_i\circ p_1\circ p_2$. \\

We will need to regularize the currents $T$ and $T_{LS}$ in a neighborhood of $0$. For that we will regularize the current $T$ and use the continuity of the Lelong-Skoda transform. But in order to do that we have to regularize the current $T$ in $X$ and not only in a neighborhood of $0$ because $T_{LS}$ is defined "globally". In the case where $X=\mathbb{P}^k$, smooth positive closed forms are dense in the space of positive closed currents. It is not true in the case of Kähler manifolds. Nevertheless, in \cite{DS5}, the authors prove the following regularization result:
\begin{theorem}\label{approx}
For a positive closed current $T$, there exist smooth positive closed forms of bidegree $(k-p,k-p)$, $T_n^+$ and $T_n^-$ such that $T_n^+-T_n^-$ converges weakly to the current $T$. Moreover, $\|T_n^{\pm}\| \leq C_X \|T\|$ where $C_X$ is a constant independent of $T$.
\end{theorem}
Let $(T_n^+)$ and $(T_n^-)$ be as in the theorem. Extracting if necessary, we suppose that the sequences converge to $T^+$ and $T^-$.  By continuity of the transform, we have $\mathcal{L}_p(T_n^{\pm})\to \mathcal{L}_p(T^{\pm})$.  Recall that the function $S\to \nu(S,0)$ is a linear form on positive closed currents. So proving the theorem for $T^{\pm}$ gives the theorem for $T$. In particular, it is enough to consider the case where $T$ is the limit in the sense of currents  of a sequence $(T_n)$ of smooth positive closed forms. We have $\|T_n\|\rightarrow \|T\|$. Then we have the lemma:
\begin{lemme}
The sequences $(\|\widehat{\pi_2}^*(T_n)\|)_n$ is bounded. In particular, we can extract a subsequence of $(\widehat{\pi_2}^*(T_n))$ that converges to a cluster value that we denote by $\widehat{T}$.
\end{lemme}
\emph{Proof.} The mass $\|\widehat{\pi_2}^*(T_n)\|$ depends only on the cohomology class $[T_n]$. The lemma follows since the cohomology class  $[T_n]$ is controlled by the mass of $T_n$ which is bounded with $n$. $\Box$ \\ 

The current $\widehat{T}$ is  positive closed as a limit of positive closed currents. As in the case of the Lelong number, the current $\widehat{T}$ is not unique. Despite that fact, we will show that the mass of $\widehat{T}$ on the set $E$ is independent of the choice of $(T_n)$ and $\widehat{T}$. More precisely, let $\widehat{T}_{|E}$ denote the restriction of $\widehat{T}$ at $E$. We have the following proposition which is the key of our proof:
\begin{proposition}\label{key}
The Lelong number $\nu(T_{LS},0)$ is equal to the mass of the current $(\Pi)_*(\widehat{T}_{|E})$ on $\Pi(E)=p_1^{-1}(0,0)$. 
\end{proposition}
 \begin{lemme}\label{3.2}
Let $S$ be a smooth positive closed form of bidimension $(p,p)$ in $X$ then $\mathcal{L}_p(S)$ is a continuous positive closed form of bidegree $(1,1)$. Furthermore:
$$ \mathcal{L}_p(S)= \widehat{\pi_1}_*(\widehat{\pi_2}^*(S)\wedge \widetilde{\Omega}^{p+1}).$$ 
\end{lemme}
\emph{Proof.} The first part of the lemma follows from the fact that $\Omega^{p+1}$ has coefficients in $L^1$. The second assertion is a consequence of the fact that $p<k-1$ hence $\Omega^{p+1}$ does not charge $\Delta$. $\Box$\\
 
\noindent We leave as an exercice to the reader the fact that $T_n$ is in fact smooth: it is a consequence of the fact that $\widetilde{\pi_1}$ and $\widetilde{\pi_1}$ are submersion. Nevertheless, continuity is sufficient for our purpose. We have the lemma:
\begin{lemme}\label{3.7} With the above notations:
$$\nu(T_{LS},0) =\lim_{r\to 0} \lim_{n \to \infty} \nu(\mathcal{L}_p(T_n),0,r).$$
\end{lemme}
\emph{Proof.}  We argue as in the proof of Lemma \ref{interpretation}: we use an approximation of the logarithm, Stokes formula and weak convergence. $\Box$\\ 

\noindent Since $\mathcal{L}_p(T_n)$ is continuous, using Lemma \ref{3.2} we have that:
$$\nu(\mathcal{L}_p(T_n),0,r)=\int_{B_r} \widehat{\pi_1}_*(\widehat{\pi_2}^*(T_n)\wedge \widetilde{\Omega}^{p+1}) \wedge (dd^c\log|x|)^{k-1}.$$
We have that $(dd^c\log|x|)^{k-1}= \widehat{\pi_1}_*(\omega_{FS, k-1}^{k-1})$ and $\widehat{\pi_1}^*((dd^c\log|x|)^{k-1})= \omega_{FS, k-1}^{k-1}$ outside $(p_1\circ p_2)^{-1}(x=0)$ (that is outside $w=0$ in the new coordinates). Let $0<r'<r$, we have that:
 \begin{eqnarray*}
\int_{B_r\backslash B_{r'}} \widehat{\pi_1}_*(\widehat{\pi_2}^*(T_n)\wedge \widetilde{\Omega}^{p+1}) \wedge (dd^c\log|x|)^{k-1}&=& \\ \int_{\widehat{\pi_1}^{-1}(B_r\backslash B_{r'})} \widehat{\pi_2}^*(T_n)\wedge \widetilde{\Omega}^{p+1}\wedge\omega_{FS, k-1}^{k-1}.
\end{eqnarray*}
So we claim that:
\begin{eqnarray}\label{formule}
\nu(\mathcal{L}_p(T_n),0,r)= \int_{\widehat{\pi_1}^{-1}(B_r)} \widehat{\pi_2}^*(T_n)\wedge \widetilde{\Omega}^{p+1}\wedge\omega_{FS, k-1}^{k-1}.
\end{eqnarray}
Indeed, the current $\widehat{\pi_1}_*(\widehat{\pi_2}^*(T_n)\wedge \widetilde{\Omega}^{p+1}) \wedge (dd^c\log|x|)^{k-1}$ does not charge $0$ and the current $\widehat{\pi_2}^*(T_n)\wedge \widetilde{\Omega}^{p+1}\wedge\omega_{FS, k-1}^{k-1}$ does not charge $w=0$.

Consider the smooth form $L$ defined by:
$$L := \widetilde{\Omega}^{p+1}\wedge \omega_{FS, k-1}^{k-1}.$$ 
As in the proof of Lemma \ref{interpretation}, letting $n\to \infty$ in (\ref{formule}) we have:
$$ \int_{\widehat{\pi_1}^{-1}(B_r)} \widehat{T}\wedge L \leq \lim_{n\to \infty} \nu(\mathcal{L}_p(T_n),0,r) \leq \int_{\widehat{\pi_1}^{-1}(B_{2r})} \widehat{T}\wedge L. $$ 
In fact we have the equality $\int_{\widehat{\pi_1}^{-1}(B_r)} \widehat{T}\wedge L = \lim_{n\to \infty} \nu(\mathcal{L}_p(T_n),0,r)$ for $r$ generic but we only need the previous inequalities. Combining this with Lemma \ref{3.7}, we have the equality:
$$ \nu(T_{LS},0)=\lim_{r\to 0} \int_{\widehat{\pi_1}^{-1}(B_r)} \widehat{T}\wedge L= \| \widehat{T}\wedge L\|_{\widehat{\pi_1}^{-1}(0)}.$$
The next lemma shows that the mass is in fact concentrated on $E$.
\begin{lemme}\label{3.8} With the above notations $\nu(T_{LS},0)$ is the mass of $\widehat{T}\wedge L$ on $E$.
\end{lemme}
\emph{Proof.} Let $W$ be a small neighborhood of $E$ in $\widehat{\widehat{X\times X}}$. It is sufficient to show that the current $\widehat{T}$ does not charge the set $\widehat{\pi_1}^{-1}(0) \backslash W$. We argue as in Lemma \ref{Fubini} taking advantage of the fact that $\widehat{\pi_2}^*$ of a current is given here by integration on fibers which are transverse to $\widehat{\pi_1}^{-1}(0)$. Of course, in $W$ the geometry is more complicated and there $\widehat{\pi_1}$ is not a submersion. $\Box$\\

\noindent \emph{End of the proof of Proposition \ref{key}.} On $\Pi(E)$, the currents $\Pi_*(L_{|E})$, $\Pi_*(\widetilde{\Omega}^{p+1}_{|E})$ and $\Pi_*((\omega_{FS,k-1}^{k-1})_{|E})$ are well defined and have no mass on $\Pi(M)$ because $M$ is of dimension $k-1$ and the bidimension of  $\Pi_*(\widetilde{\Omega}^{p+1}_{|E})$ and $\Pi_*((\omega_{FS,k-1}^{k-1})_{|E})$ is at least $k$ and the singularities are in two disjoint subvarieties ($\Pi_*(\widetilde{\Omega}^{p+1}_{|E})$ is smooth where $\Pi_*((\omega_{FS,k-1}^{k-1})_{|E})$ is singular and vice versa). Furthermore, $ \Pi_*(L_{|E})=\Pi_*(\widetilde{\Omega}^{p+1}_{|E}) \wedge \Pi_*((\omega_{FS,k-1}^{k-1})_{|E})$.  On $\Pi(E)$, the measure $\Pi_*(\widehat{T}_{|E}\wedge L_{|E})$ is well defined as its trivial extension over $E$ of its restriction to $E\backslash M$. Indeed, $\widehat{T}_{|E}$ does not charge the set $(p_2)^{-1}(M)$ (it is the same argument as in Lemma \ref{Fubini}). In particular, the measure is equal to $\Pi_*(\widehat{T}_{|E})\wedge \Pi_*( L_{|E})$ and, by Lemma \ref{3.8}, its mass is equal to $\nu(T_{LS},0)$.  

 Then the mass of the measure can be computed in cohomology since for positive closed currents in $\mathbb{P}^{2k-1}$ the mass of a wedge product is the product of the masses. Let $F$ be a subspace of $\mathbb{C}^{2k}$ and consider the orthogonal projection from $\mathbb{C}^{2k}$ to $F$. It induces a meromorphic map  $\sigma_F:\mathbb{P}^{2k-1} \dashrightarrow \mathbb{P}(F)$. Then for a positive closed current $S$ of mass $m$ on $\mathbb{P}(F)$, the pull-back $\sigma_F^*(S)$ is well defined and of mass $m$ (see Section 1 in \cite{Meo}). Applying this to $\mathbb{P}(F_1)=[v=0]$ and $\mathbb{P}(F_2)=[u=0]$ we obtain that $\Pi_*( L_{|E})= \Pi_*( \widetilde{\Omega}^{p+1}_{|E})\wedge\Pi_*( (\omega^{k-1}_{FS,k-1})_{|E})$ is of mass $1$. Indeed $\widetilde{\Omega}^{p+1}_{|E}$ is by definition $\widetilde{\Omega}^{p+1}_{|E}(0,0,[u'])$. So
\begin{eqnarray*}
\Pi_*( \widetilde{\Omega}^{p+1}_{|E})&=& \sigma_{F_1}^*(\widetilde{\Omega}^{p+1}), 
\end{eqnarray*}
since the equality is true outside $\mathbb{P}(F_2)$ and both sides of the equality give no mass to $\mathbb{P}(F_2)$. Similarly:
\begin{eqnarray*}
\Pi_*( (\omega_{FS,k-1}^{k-1})_{|E})   &=& \sigma_{F_2}^*(\omega_{FS,k-1}^{k-1}). 
\end{eqnarray*}
And our normalization of $\widetilde{\Omega}$ implies that $\Pi_*( \widetilde{\Omega}^{p+1}_{|E})$ is of mass $1$ so $\Pi_*( L_{|E})$ is of mass $1$. Then $\nu(T_{LS},0)$ can be interpreted as the mass of $\Pi(E)$ for the current $\Pi_*(\widehat{T}_{|E})$ which is also the mass of $(p_2)_*(\widehat{T})$ on $\Pi(E)$. $\Box$\\

\noindent \emph{End of the proof of Theorem \ref{main}.}  The Lelong number $\nu(T,0)$ is the same than the Lelong number $\nu(\pi_2^*(T),(0,0))$ (more generally, the Lelong numbers of the pull-back of a current by a submersion are preserved by Proposition (2.3) in \cite{Meo}). This and Lemma \ref{interpretation} applied to $\widetilde{T}:= (p_2)_*(\widehat{T})=\lim_n p_1^*(T_n) $ imply that $\nu(T,0)$ is the mass of $\Pi(E)$ for the current $(p_2)_*(\widehat{T})$ (we use that $(p_2)_* (p_2)^*= \text{id}$). Proposition \ref{key} implies the result. $\Box$\\

\noindent Several remarks are in order here.
\begin{Remark} \rm
The transform $\mathcal{L}_p$ is compatible with the cohomology, that is if $T$ and $T'$ are cohomologuous on $H^{p,p}(X,\mathbb{C})$ so are $\mathcal{L}_p(T)$ and $\mathcal{L}_p(T')$ by Lemma \ref{Fubini}. Indeed, if $T=T'+dd^c \alpha$ with $\alpha$ of bidegree $(k-p-1, k-p-1)$ then since the $\widetilde{\pi_i}$ are submersions and $\widetilde{\Omega}$ is closed we have:
$$\mathcal{L}_p(T)-\mathcal{L}_p(T')= dd^c (\widetilde{\pi_1}_*(\widetilde{\pi_2}^*(\alpha) \wedge \widetilde{\Omega}^{p+1})).$$  
 \end{Remark}
\begin{Remark} \rm
The choice $\widetilde{\Omega}^{p+1}$ can be replaced by a strongly positive closed smooth form $\Theta$ on $\widetilde{X\times X}$ of bidegree $(p+1,p+1)$ such that the mass of $\Theta\wedge [\pi_2^{-1}(y)]$ is equal to $1$ for any $y$ (this mass is a constant for cohomological reasons, so we just have to normalize it). 
\end{Remark}
\begin{Remark} \rm
The same method allows us to prove that: "for any current $T$ of bidimension ($p,p)$ ($p<k-1$) and any $p<q \leq k-1$, there exists a positive closed current $T_q$ of bidimension $(q,q)$ depending continuously and linearly of $T$ which has the same Lelong number as $T$ at every point".
\end{Remark}

\section{Generalized Demailly's inequality}
Using Theorem \ref{main}, we follow the argument of Méo (\cite{Meo}) to prove Theorem \ref{1}. We use the notations of the introduction. Demailly proved the following regularization result:
\begin{theorem}[\cite{dem1}]
Let $S$ be a positive closed current of bidegree $(1,1)$ on a compact Kähler manifold $X$. Let $\varphi$ be a quasi-psh function such that $S=\alpha+dd^c\varphi$ where $\alpha$ is a smooth $(1,1)$ form. Then for all $c>0$, there exists a decreasing sequence $(\varphi_{c,l})_{l\geq 1}$ of functions converging to $\varphi$ such that:
\begin{itemize}
\item $\varphi_{c,l}$ is smooth outside $X \backslash E_c$ ;
\item $dd^c \varphi_{c,l}+\alpha +A \|S\| \omega \geq 0$ where $A$ is a constant that depends only on $X$ and $\omega$;
\item for all $x\in X$, $\nu(\varphi_{c,l}, x)= ( \nu(\varphi, x)-c)_+ := \max( \nu(\varphi, x)-c, 0)$.
\end{itemize}
\end{theorem}
We apply that result for $S=T_{LS}$. The current:
 $$ T \wedge (dd^c \varphi_{c_{p-1},l}+\alpha+A \|S\| \omega) \wedge \dots \wedge (dd^c \varphi_{c_q,l}+\alpha+A \|S\| \omega) $$
is well defined by the theory of intersection of currents (see \cite{dem3} and \cite{FS2}) because for $c_j > b_j$, the set of points at which $\varphi_{c_j,l}$ is not bounded is contained in $E_{c_j,l}$ which is of dimension less or equal to $j$. Using \cite[Corollary (7.9) p. 194]{dem2}, its Lelong number at $x$ is greater or equal than:
$$ \nu(T,x)(\nu(T,x)-c_{p-1})_+\dots(\nu(T,x)-c_{q})_+.$$
Siu's theorem \cite{Siu} implies that this current is greater than:
 $$\sum_n \nu_{q,n}(\nu_{q,n}-c_{p-1})_+\dots(\nu_{q,n}-c_{q})_+ [Z_{q,n}].$$
Observe that the mass of the current $dd^c \varphi_{c,l}+\alpha +A \|T\| \omega \geq 0$ is equal to $(1+A)\| T\|$. So taking the masses gives:
$$ \sum_n \nu_{q,n}(\nu_{q,n}-c_{p-1})_+\dots(\nu_{q,k}-c_{q})_+ \|[Z_{q,n}]\| \leq C \|T\|^{p-q+1}.$$
Theorem \ref{1} follows from letting the $c_j$ go to $b_j$.

\section{Transformation of pluriharmonic currents}
We want to generalize the results of Sections \ref{geometric} and \ref{properties} to the case of pluriharmonic currents. Once again we write "linear operator" on a convex cone instead of an affine operator on a convex cone which extends to a linear operator on the vector space it spans. 
\begin{theorem}
The Lelong-Skoda transform $\mathcal{L}_p$ is a well defined linear continuous operator from the cone of positive pluriharmonic currents of bidimension $(p,p)$ ($p<k-1$) to the cone of positive pluriharmonic currents of bididegree $(1,1)$. The transform preserves Lelong numbers.
\end{theorem}
Recall that a positive current $T$ of bidimension $(p,p)$ is said to be \emph{pluriharmonic} if $dd^c T=0$ in the sense of currents (see \cite{FS3}). For such current, Skoda proved that the Lelong number is well defined (Proposition 1 in \cite{Sko1}). More precisely, let $T$ be a positive pluriharmonic current of bidimension $(p,p)$ in an open set $U$ of $\mathbb{C}^k$. Then, for $x$ in $U$, we have that the positive measure $T(y) \wedge (d_yd_y^c  \log |x-y|)^p$ is well defined  on $U \backslash \{x\}$ since $\log |x-y|$ is smooth here. And for $0<r_1<r_2$ we have the identity:
\begin{eqnarray}\label{Skoda}
  \frac{1}{r_2^p}\int_{|x-y|<r_2} T(y)\wedge (\frac{1}{2}dd^c |y|^2)^p-\frac{1}{r_1^p}\int_{|x-y|<r_1} T(y)\wedge (\frac{1}{2}dd^c |y|^2)^p=  \nonumber\\
\int_{r_1<|x-y|<r_2} T(y) \wedge (d_yd_y^c  \log |x-y|)^p . 
\end{eqnarray}
In particular, the non negative quantity $\frac{1}{r^p}\int_{|x-y|<r} T(y)\wedge (\frac{1}{2}dd^c |y|^2)^p$ decreases with $r$ to a number $\nu(T,x)$ called the \emph{Lelong number} of $T$ at $x$. Let $(T_n)$ be a sequence of smooth positive pluriharmonic currents converging to $T$ in the sense of currents in a neighborhood of $x$. Then since $\lim_{n\to  \infty} T_n(y)\wedge (\frac{1}{2}dd^c |y|^2)^p \to T(y)\wedge (\frac{1}{2}dd^c |y|^2)^p$ in the sense of measures, we have for $r$ generic that:
$$\lim_{n\to \infty} \frac{1}{r^p}\int_{|x-y|<r} T_n(y)\wedge (\frac{1}{2}dd^c |y|^2)^p \to \frac{1}{r^p}\int_{|x-y|<r} T(y)\wedge (\frac{1}{2}dd^c |y|^2)^p. $$   
For $T_n$ smooth, we can let $r_1$ goes to $0$ in Skoda's formula (\ref{Skoda}) and we see that:   
$$\frac{1}{r^p}\int_{|x-y|<r} T_n(y)\wedge (\frac{1}{2}dd^c |y|^2)^p= \int_{|x-y|<r} T_n(y) \wedge (d_yd_y^c  \log |x-y|)^p. $$
Combining the last two equalities gives: 
\begin{eqnarray}\label{harmo}
 \nu(T,x)= \lim_{r\to 0} \lim_{n \to \infty} \int_{|x-y|<r} T_n(y) \wedge (d_yd_y^c  \log |x-y|)^p . 
 \end{eqnarray}
In the case of pluriharmonic currents, we have the following integration by parts formula (see \cite{dem3}) :
\begin{lemme} Let $T$ be a pluriharmonic current of bidimension $(p,p)$ on $U \subset X$ and $f$ be a smooth form of bidegree $(p-1,p-1)$ equal to zero near $\partial U$, then:
$$\int_U T \wedge dd^c f=0.$$  
\end{lemme}
Once again, the Lelong number at a point does not depend on the choice of coordinates (see \cite{AB1} for the proof in a more general setting of pluriharmonic currents), so we can speak of Lelong numbers on a manifold. Lemma \ref{interpretation} still applies for positive pluriharmonic currents. That is: if $\imath:\widehat{X}\to X$ is the blow-up of $X$ at $x$ and $\widehat{T}$ is a cluster value of the sequence $\imath^*(T_n)$, then the mass of $\widehat{T}$ on the exceptional divisor is equal to the Lelong number $\nu(T,x)$ (the proof is exactly the same using formula (\ref{harmo})). Now we define in the same way $\mathcal{L}_p(T)$ which is a well defined positive pluriharmonic current of bidegree $(1,1)$. Observe that the arguments of Section \ref{properties} remain valid so we can conclude in the same way.

 Let us mention the points where the argument need some modifications. In \cite{DS5}, the authors also prove Theorem \ref{approx} for positive pluriharmonic currents. In Lemma \ref{3.7}, we use the previous argument instead of Stokes formula: we use formula (\ref{harmo}) applied to $\mathcal{L}_p(T_n)$ instead of formula (\ref{trick}) ($\mathcal{L}_p(T_n)$ is continuous and not smooth but that is enough here). To prove that the measures $T_{|E}\wedge L_{|E}$ and $T_{|E}\wedge \omega^{k+p}_{FS, 2k-1}$ have the same mass, observe that $\omega^{k+p}_{FS, 2k-1}$ and  $L_{|E}$ are cohomologous since they have the same mass so $\omega^{k+p}_{FS, 2k-1}-L_{|E}=dd^c f$ where $f$ is a form with coefficient in $L^1$. And we conclude using the previous integration by part formula.\\

\noindent In an open set $U$ of $\mathbb{C}^k$, it is not true that the set $E_c=\left\{ z\in U, \ \nu(T,y)>c \right\}$ is analytic for a pluriharmonic current $T$. Indeed, consider in $\mathbb{C}^2$  the current $T= h(z_1) [z_2=0]$ where $h$ is a non-constant non-negative harmonic function in $\mathbb{C}$ and $[z_2=0]$ is the current of integration on $\{z_2=0\}$. In a compact manifold, any pluriharmonic function is constant. That raises the question:\\

\noindent {\bf Open problem} (Dinh): Let $T$ be a positive pluriharmonic current on a compact Kähler manifold $(X,\omega)$. Does Siu's theorem still hold for $T$, i.e. is $E_c$ analytic for $c>0$?\\

A relative result was proved in \cite{DL1} in the case of rectifiable currents. The previous theorem simplifies the question to the case of bidegree $(1,1)$.

 \noindent Gabriel Vigny, Mathématiques - Bât. 425, UMR 8628,\\ 
 Université Paris-Sud, 91405 Orsay, France.  \\
\noindent Email: gabriel.vigny@math.u-psud.fr

\end{document}